\documentclass[11pt]{article}
\usepackage{amsthm, amsmath, amssymb}

\newtheorem{prop}{Proposition}
\newtheorem{thm}{Theorem}

\newtheorem{lemma}{Lemma}
\newtheorem{conj}{Conjecture}

\title{The characteristic ideal of a \\ finite, connected, regular graph
\thanks{Work partially supported by Ministerio de Ciencia y
Tecnolog\'{\i}a under projects BFM2001-2340 and BFM 2003-00368 and by
Generalitat de Catalunya under project 2001 SGR 00224.}}

\author{
Josep M. Brunat \& Antonio Montes
       \\ Departament de Matem\`atica Aplicada 2 \\
       Universitat Polit\`ecnica de Catalunya, Spain\\
       \{Josep.M.Brunat,Antonio.Montes\}@upc.edu\\
       http://www-ma2.upc.edu/$\sim$montes}
       \date{July, 2004}
\begin{document}
%
\setcounter{page}{1}

\maketitle
\begin{abstract}
  Let $\Phi(x,y)\in\mathbb{C}[x,y]$ be a symmetric polynomial of partial
  degree $d$. The
  graph $G(\Phi)$ is defined by taking $\mathbb{C}$ as set of vertices
  and the points of $\mathbb{V}(\Phi(x,y))$ as edges. We study the following
  problem: given a finite, connected, $d$-regular graph $H$, find the
  polynomials $\Phi(x,y)$ such that $G(\Phi)$ has some connected
  component isomorphic to $H$ and, in this case, if $G(\Phi)$ has (almost)
  all components isomorphic to $H$. The problem is solved by associating to
  $H$ a characteristic ideal which offers a new perspective to the
  conjecture formulated in a previous paper, and allows to reduce its
  scope.  In the second part, we determine the characteristic ideal
  for cycles of lengths $\le 5$ and for complete graphs of order $\le 6$.
  This results provide new evidence for the conjecture.
\end{abstract}




\section{Introduction}
In the previous papers~\cite{BrMo,BrLa} there are given
basic notations and descriptions that will be assumed here and we
refer to them for all not defined concepts. In this paper, we only
consider symmetrical polynomials. Let us recall two basic
definitions restricted to a symmetrical polynomial
$\Phi(x,y)\in\mathbb{C}[x,y]$ of partial degree $d$. The graph
$G(\Phi)$ is defined by taking $\mathbb{C}$ as set of vertices and
the points of $\mathbb{V}(\Phi(x,y))$ as edges. As shown
in~\cite{BrMo}, for standard symmetrical polynomials of partial degree $d$
(defined in~\cite{BrMo}), all the connected components of $G(\Phi)$
but a finite number are $d$-regular graphs without loops nor
multiple arcs nor defective vertices. The graph $G(\Phi)^*$ is obtained by
removing from $G(\Phi)$ the finite set of singular components.

The problem studied here is the following: given a finite,
connected, $d$-regular graph $H$, find the polynomials $\Phi(x,y)$
(if any exists),  such that $H$ is isomorphic to some (connected)
component of $G(\Phi)^*$. If it is the case, the question of
deciding when $H$ is isomorphic to \emph{all} components of
$G(\Phi)^*$ is the matter of the conjecture formulated
in~\cite{BrMo}, which, for symmetric polynomials, admits the
following formulation: If $\Phi(x,y)\in\mathbb{C}[x,y]$ is a
standard symmetric polynomial and $G(\Phi)^*$ has a finite
component, then all components are isomorphic.

In Section 2, we define a system $S(H,\Phi)$ and a variety
$W(H,\Phi)$ associated to a \emph{pairing}, which is a pair
$(H,\Phi)$ formed by a finite, connected, $d$-regular graph $H$,
and a symmetric polynomial of partial degree $d$. The points
$(u_1,\ldots,u_n)$ of $W(H,\Phi)$ such that $u_1,\ldots,u_n$
induce a component in $G(\Phi)$ form a variety $U(H,\Phi)\subseteq
W(H,\Phi)$. The correspondence between points of $U(H,\Phi)$ and
components of $G(\Phi)$ is established.

In Section 3  we characterize the finite, connected, $d$-regular
graphs that are isomorphic to a component of some $G(\Phi)^*$ by
its associated {\em characteristic ideal}. This leads to an
algebraic formulation of the conjecture, and to the reduction of
its scope.  It also provides the theoretic frame for constructing
an algorithm to determine the characteristic ideal of $H$.

In Section 4 we show how to improve the initial polynomial system
by eliminating undesired solutions in order to determine the
characteristic ideal of a graph. The general algorithm is applied
to find the characteristic ideals of $3$-cycles and $4$-cycles.

Because of the complexity of the computations using the general
algorithm, specific algorithms are valuable for some kind of graphs.
In section~5, we give an algorithm for cycles, providing the
characteristic ideal for cycles of length $\le 5$; and, in Section~6
another for complete graphs, providing the characteristic ideal for
complete graphs of order $\le 6$. All these results provide further
evidence of the conjecture, besides those obtained in~\cite{BrMo}.

Finally in the conclusions, some open problems are formulated.

Besides~\cite{BrMo}, for undefined algebraic concepts we refer
to~\cite{CoLiSh,CoLiSh2}, and for graph theoretic ones to~\cite{ChLe,Watkins2}.

\section{The variety of a pairing}
\label{AlgCharacterization}

Let $H$ be a finite, connected, $d$-regular graph and $\Phi(x,y)$
a symmetric polynomial of partial degree $d$. The immediate goal
is to decide if $H$ is isomorphic to a component of $G(\Phi)^*$. A
\emph{pairing} (or \emph{$d$-pairing} if we wish emphasize the
degree of the graph and the partial degree of the polynomial) is a
pair $(H,\Phi)$ where $H$ is a finite, connected, $d$-regular
graph ($d\ge 1$) and $\Phi(x,y)$ a symmetric polynomial of partial
degree $d$. A pairing $(H,\Phi)$ is \emph{standard} (resp.
\emph{non
  standard}) if the polynomial $\Phi(x,y)$ is standard (resp. non
standard). We shall assume that, if $H$ is of order $n$, the set of
vertices of $H$ is $[n]=\{1,\ldots,n\}$.  Associate to a $d$-pairing
$(H,\Phi)$, with $H=([n],E)$, we define the set of polynomials
$$
S=S(H,\Phi)=\{\Phi(x_i,x_j)\ :\ ij\in E\},
$$
and the variety of $S(H,\Phi)$,
$$
W=W(H,\Phi)=\mathbb{V}(S(H,\Phi)).
$$
Note that if $H$ is $d$-regular of order $n$, then $H$ has $m=dn/2$
edges, and the system $S$ has $m$ polynomials.
Moreover $W\subseteq \mathbb{C}^n$.
A point $(u_1,\ldots,u_n)\in W$ is
called a \emph{proper point} if $u_i\ne u_j$ for $1\le i<j\le n$;
otherwise it is an \emph{improper} point.

\begin{lemma}
\label{dimW}
Let $(H,\Phi)$ be a pairing. Then, the dimension of the variety
$W(H,\Phi)$ is at most 1.
\end{lemma}
\begin{proof}
  Let $H$ be of order $n$ and degree $d$. Let $(u_1,\ldots,u_n)$ be a
  point in $W$. If $1j$ is an edge of $H$, then the number of
  distinct values of $x_j$ in the points of $W$ with $x_1=u_1$
  is at most $d$, the maximum number of roots of $\Phi(u_1,y)$. By
  induction, if $j$ is a vertex of $H$ at distance $r$ from the vertex
  1, then the number of distinct values of $x_j$ in points of
  $W$ with $x_1=u_1$ is at most $d^r$. Therefore, the number of
  points of $W$ with $x_1=u_1$ is finite.
\end{proof}

The following theorem shows that for standard pairings $(H,\Phi)$,
the proper points of $W(H,\Phi)$ correspond to components of $G(\Phi)^*$.
Given $u_1,\ldots,u_n\in\mathbb{C}$, we denote by $\langle u_1,\ldots
u_n \rangle$ the subdigraph of $G(\Phi)$ induced by $u_1,\ldots,u_n$
(the polynomial $\Phi(x,y)$ being implicit).

\begin{prop}
\label{proper}
Let $(H,\Phi)$ be a standard pairing. Then
a point $(u_1,\ldots,u_n)\in W(H,\Phi)$ is a proper point if and
only if $\langle u_1,\ldots u_n \rangle$ is a
component of $G(\Phi)^*$ isomorphic to $H$.
\end{prop}
\begin{proof}
  Let $(u_1,\ldots,u_n)\in W$ be a proper point.
  Define $f\colon [n]\rightarrow \{u_1,\ldots,u_n\}$
  by $f(i)=u_i$. For $i\ne j$ we have $u_i\ne
  u_j$, hence $f$ is injective.  As the two sets $[n]$ and
  $\{u_1,\ldots,u_n\}$ have the same cardinality $n$, the mapping $f$
  is bijective.

  If $ij\in E$, then $(u_i,u_j)$ is a zero of the polynomial in
  $S$ corresponding to the edge $ij$, that is, $u_i$ is
  adjacent to $u_j$ in $G(\Phi)$. As $H$ is connected, the subdigraph
  $\langle u_1,\ldots,u_n\rangle$ is connected.
  Both graphs are $d$-regular, so $f$ is an
  isomorphism.  From the fact that $H$ is a $d$-regular graph, it
  follows that it has neither loops, nor
  multiple edges, nor defective vertices. Therefore $\langle
  u_1,\ldots,u_n\rangle=G(\Phi, u_1)$ is not a singular component of $G(\Phi)$.

  Conversely, let $C$ be a component of $G(\Phi)^*$ isomorphic to $H$.
  Let $f\colon i\mapsto u_i$ be the isomorphism from $H$ onto $C$.
  Clearly, if $1\le i<j\le n$, then $u_i\ne u_j$. For each polynomial
  $\Phi(x_i,x_j)$ of $S$, we have an edge $ij\in E$. As $f$ is
  an isomorphism, $u_i$ is adjacent to $u_j$ in $G(\Phi)$, which is
  equivalent to $\Phi(u_i,u_j)=0$. Therefore, $(u_1,\ldots,u_n)$ is a
  proper solution of $S(H,\Phi)$.
  \end{proof}

  Now consider improper points of $W$. Recall that, even if
  $\Phi(x,y)$ is a symmetric polynomial, the singular components of
  $G(\Phi)$ can be digraphs with loops or multiple arcs. The following
  decomposition helps to eliminate solutions of $S$ which do
  not correspond to components of $G(\Phi)$. For a given pairing
  $(H,\Phi)$, define
$$
\begin{array}{rcl}
Z=Z(H,\Phi) &=&  W(H,\Phi)
               \bigcap
               \left( \bigcup_{i>j} \mathbb{V} (x_i-x_j) \right); \\
U=U(H,\Phi) &=& \overline{ W(H,\Phi)
              \setminus
              Z(H,\Phi)}.
\end{array}
$$
Note that $Z$ is the set of improper points of $W$, and that the
proper points of $W$ are in $W\setminus Z$, so they are in its
algebraic closure $U$.

\begin{prop}
\label{J}
Let $(H,\Phi)$ be a standard $d$-pairing and let $J$ be the set of
improper points of $U(H,\Phi)$. Then
\begin{enumerate}
\item[\emph{(i)}] The set $J$ is finite.
\item[\emph{(ii)}] If $J\ne\emptyset$, then $\operatorname{dim} U(H,\Phi)=1$.
\end{enumerate}
\end{prop}
\begin{proof}
(i) If $J$ is empty, the result is trivial. Assume $J\ne\emptyset$.
Note that $J=U\cap\left( \bigcup_{j<i} \mathbb{V} (x_i-x_j) \right)$.
We have $J\subseteq U\subseteq W$ and, by Lemma~\ref{dimW},
$\operatorname{dim}W\le 1$. Therefore, $\operatorname{dim} J\le 1$.
If $J$ were not finite, then
$1=\operatorname{dim}J\le\operatorname{dim}U\le 1$. If $U$ is
irreducible, then $J=U$. In this case, $W\setminus Z\subseteq U=J$,
and the set $W\setminus Z$ has no proper points, so it is empty.
Therefore $J=U=\overline{W\setminus
  Z}=\overline{\emptyset}=\emptyset$, a contradiction. If $U$ is reducible,
then decompose it as a union of irreducible varieties $U=\bigcup_{i=1}^s U_i$
and set $J_i=J\cap U_i$. For each infinite $J_i$, we have $J_i=U_i$
and $\overline{U_i\setminus Z}=\emptyset$. Therefore
$U=\overline{W\setminus Z}=\bigcup_{i=1}^s \overline{U_i\setminus Z}$
is finite, a contradiction.

(ii) Let $\overline{u}\in J$. We have $\overline{u}\notin
W\setminus Z$ but $\overline{u}\in \overline{W\setminus Z}$. Thus
$\overline{u}$ is not an isolated point of $U$. Therefore $U$ is
infinite.  By using Lemma~\ref{dimW}, we have $1\le
\operatorname{dim}U\le \operatorname{dim}W \le 1$. Therefore, we
conclude $\operatorname{dim}U=1$.
\end{proof}

As a consequence of propositions~\ref{proper} and \ref{J}, we have:
\begin{thm}
\label{HU}
Let $(H,\Phi)$ be a standard pairing. Then the graph $G(\Phi)^*$ has
some component isomorphic to $H$ if and only if $U(H,\Phi)\ne\emptyset$.
\end{thm}
\begin{proof}
Proposition~\ref{proper} ensures that if $G(\Phi)^*$ has a component
$C$ isomorphic to $H$ and $i\mapsto u_i$ is the isomorphism from $H$ to
$C$, then $(u_1,\ldots,u_n)$ is a proper solution, i.e.,
$(u_1,\ldots,u_n)\in W\setminus Z\subseteq U$. Therefore $U\ne
\emptyset$.

Conversely, assume that there exists
$\overline{u}=(u_1,\ldots,u_n)\in U$. If $\overline{u}$ is a
proper point, then Proposition~\ref{proper} ensures that a
component of $G(\Phi)^*$ is isomorphic to $H$. If $\overline{u}$
is an improper point, then by Proposition~\ref{J},
$\operatorname{dim}U=1$. Therefore, $W\setminus Z$ is not empty
and there exists a proper point $\overline{u}\in W\setminus Z$. By
Proposition~\ref{proper}, there exists a component of $G(\Phi)^*$
isomorphic to $H$.
\end{proof}

Consider now non standard pairings.

\begin{prop}
\label{HUns} Let $(H,\Phi)$ be a non standard $d$-pairing. Put
$n=d+1$. If $U(H,\Phi)\ne\emptyset$ then
\begin{enumerate}
\item [\emph{(i)}] $\Phi(x,y)=f(x)f(y)\Phi_1(x,y)$ where $\operatorname{deg}f(x)\ge
  1$ and $\Phi_1(x,y)$ is a standard polynomial.
\item [\emph{(ii)}]$G(\Phi)$ has universal vertices, say
  $w_1,\ldots,w_r$, and it is connected.
\item [\emph{(iii)}] All non singular components of
  $G(\Phi_1)^*$ are isomorphic to $K_{n-r}$,
  and $H$ is isomorphic to $K_{n}$.
\end{enumerate}
\end{prop}
\begin{proof}
  (i) If $L(x)=\Phi(x,x)$ is the zero polynomial or
  $\Phi(x,y)\ne\operatorname{rad}\Phi(x,y)$ then each point
  $\overline{u}\in W$ has some repeated coordinates. Hence,
  $\overline{u}\in Z$. Then, $W\setminus
  Z=\emptyset$ and $U=\emptyset$. As $\Phi(x,y)$ is
  non standard, it must be of the form $\Phi(x,y)=f(x)f(y)\Phi_1(x,y)$
  with $\operatorname{deg}f(x)\ge 1$ and $\Phi_1(x,y)$ standard.

  (ii) The roots $w_1,\ldots,w_r$ of $f(x)$ are the universal vertices.
  The existence of universal vertices implies that $G(\Phi)$ is connected.

  (iii) A point $(u_1,\ldots,u_n)\in U$ has $r$ coordinates which are
  the $r$ universal vertices. The remaining coordinates induce
  a subgraph $(n-1-r)$-regular. The partial degree of $\Phi_1(x,y)$ is
  $d-r$. Therefore $n-1=d$ and the components of $G(\Phi_1)^*$ are
  isomorphic to $K_{n-r}$. Moreover $H\simeq \langle
  u_1,\ldots,u_n\rangle=K_{n}$.
\end{proof}

In the context of graphs, the conjecture stated in~\cite{BrMo} is the
following:
\begin{conj}
\label{ccc}
Let $(H,\Phi)$ be a standard pairing. If $H$ is
isomorphic to a component of $G(\Phi)^*$, then $H$ is isomorphic to
all components of $G(\Phi)^*$.
\end{conj}

Let $(H,\Phi)$ be a standard pairing. The graph $H$ is said to be
$\Phi$-\emph{polynomial} if it is isomorphic to a component of
$G(\Phi)^*$; $H$ is said to be \emph{strongly $\Phi$-polynomial} if it
is isomorphic to all components of $G(\Phi)^*$. Conjecture~\ref{ccc}
states that if $H$ is $\Phi$-polynomial, then $H$ is strongly
$\Phi$-polynomial.

A finite, connected, $d$-regular graph $H$ is
\emph{polynomial} (resp. \emph{strongly polynomial}) if it is
$\Phi$-polynomial (resp. \emph{strongly $\Phi$-polynomial}) for some
standard polynomial $\Phi(x,y)$.

The condition of being strongly polynomial graph is quite
restrictive. Indeed, only vertex-transitive graphs can be strongly
polynomial, as shown in the following theorem.

\begin{thm}
\label{vt}
Let $H$ be a strongly $\Phi$-polynomial graph.  Then $G(\Phi)^*$ is
vertex-transitive. In particular, $H$ is vertex transitive.
\end{thm}
\begin{proof}
  Each component of $G(\Phi)^*$ is isomorphic to $H$ and, by
  Proposition~\ref{proper}, each component provides a proper point of
  $W$. The number of components of $G(\Phi)^*$ is uncountable so, by
  Lemma~\ref{dimW}, $\operatorname{dim}W=1$. Therefore, in the
  system $S$, one indeterminate, say $x_1$ is free.  For each vertex
  $u_1$ of $G(\Phi)^*$, we have some proper point of $W$ of the form
  $(u_1,\ldots,u_n)$ and an isomorphism $f_{u_1}$ from $H$ to $G(\Phi,
  u_1)$ given by $i\mapsto u_i$.  Let $u_1,v_1\in\mathbb{C}$ be two
  vertices in $G(\Phi)^*$. Then $f=f_{v_1}f_{u_1}^{-1}$ is an
  isomorphism from $G(\Phi,u_1)$ onto $G(\Phi,v_1)$ which applies
  $u_1$ in $v_1$. This implies that $G(\Phi)^*$ is vertex transitive.
  In particular, each component of $G(\Phi)^*$, which is isomorphic to
  $H$, is vertex transitive.
\end{proof}

Thus, only finite, connected, $d$-regular, vertex symmetric graphs
can be strongly polynomial. On the other side we cannot ensure
that every finite, connected, $d$-regular, vertex symmetrical
graph is strongly polynomial. Petersen's graph is the smallest
$d$-regular vertex transitive graph for which we do not know if it
is polynomial.  Our guess is that it is not, but the question is
not yet settled. All the strongly polynomial graphs given
in~\cite{BrMo} are Cayley graphs. The fact that Petersen's graph
is a well-known example of a vertex transitive graph which is not
a Cayley graph suggests that it is possible that every strongly
polynomial graph is not only vertex transitive, as
Theorem~\ref{vt} ensures, but also a Cayley graph.

\section{The characteristic ideal of a graph}

Fix a finite, connected, $d$-regular, graph $H=([n],E)$. If the
goal is to find polynomials $\Phi(x,y)$ such that $H$ is
isomorphic to one or all components of $G(\Phi)^*$, the
coefficients of $\Phi(x,y)$ must be unknowns. Then we define
$S(H)$, $W(H)$, $Z(H)$ and $U(H)$ in a similar way as in the
previous section, but considering the coefficients of the
polynomials also as unknowns. Let $m=(d+1)(d+2)/2$. For each
$\overline{a}=
(a_{d\,d},a_{d\,d-1},\ldots,a_{d\,0},a_{d-1\,d-1},\ldots,a_{0\,0})
\in\mathbb{C}^m$ let
$$
\Phi_{\overline{a}}(x,y)=\sum_{i,j=1}^d a_{i\,j}x^iy^j, \hbox{\ where\ } a_{i\,j}=a_{j\,i}.
$$
As before, define
$$
\begin{array}{cll}
S(H) &=& \{\Phi_{\overline{a}}(x_i,x_j) \ : \ ij\in E\},\\
W(H) &=& \mathbb{V}(S(H))\subseteq\mathbb{C}^{m+n}, \\
Z(H) &=& W(H)
          \bigcap
          \left( \bigcup_{i>j} \mathbb{V} (x_i-x_j) \right), \\
U(H) &=& \overline{W(H)\setminus Z(H)}.
\end{array}
$$
A point $(\overline{c},\overline{u})=(c_{d\,d},\ldots,c_{0\,0},u_1,\ldots,u_n)$
of $W(H)$ is said to be a
\emph{proper point} if  $c_{d\,j}\ne 0$ for some $j$ and $\overline{u}$
is a proper point of $S(H,\Phi_{\overline{c}})$; otherwise
it is an \emph{improper point}.

In order to decide wether $H$ is polynomial or not, the following
ideals are the key. Define
$$
\begin{array}{lll}
\mathcal{I}(H) &=& \mathbb{I}(U(H)), \\
\mathcal{I}_{\overline{a}}(H)
   &=&
\mathcal{I}(H)\cap\mathbb{C}[\overline{a}], \\
\mathcal{I}_{\overline{a},x_1}(H)
   &=&
\mathcal{I}(H)\cap\mathbb{C}[\overline{a},x_1].
\end{array}
$$
These three ideals satisfy $\mathcal{I}_{\overline{a}}(H)
\subseteq \mathcal{I}_{\overline{a},x_1}(H) \subseteq
\mathcal{I}(H)$. The ideal $\mathcal{I}_{\overline{a}}(H)$ is
called the \emph{characteristic ideal} of $H$, its name being
justified by Theorem~\ref{characideal}. First, let us put aside a
special case.

If $H=K_{n}$, then $H$ is circulant, hence strongly polynomial,
see~\cite{BrMo}. On the other hand, Proposition~\ref{HUns} shows
that there exists a non standard polynomial $\Phi(x,y)$ and a
point $(u_1,\ldots,u_n)\in U(H,\Phi)$ such that $G(\Phi)$ is
connected and $\langle u_1,\ldots,u_n\rangle$ is isomorphic to
$H$. Thus, we may consider only the case $H\ne K_{n}$.

\begin{thm}
\label{characideal} Let $H$, $(H\ne K_{d+1})$, be a finite,
connected, $d$-regular graph. Then, one of the three following
statements
holds. 
\begin{enumerate}
\item[\emph{(i)}] $\mathcal{I}(H)=\mathcal{I}_{\overline{a}(H)}=\langle 1\rangle$.
In this case $H$ is not a polynomial graph.
\item[\emph{(ii)}] $\mathcal{I}(H)\ne \langle 1\rangle$
      and
      $\mathcal{I}_{\overline{a},x_1}(H)=\mathcal{I}_{\overline{a}}(H)$.
  In this case, for all $\overline{c}\in\mathbb{V}(I_{\overline{a}})$
  the polynomial $\Phi_{\overline{c}}(x,y)$ is standard, and $H$ is a
  strongly $\Phi_{\overline{c}}$-polynomial graph.
\item[\emph{(iii)}] $\mathcal{I}(H)\ne \langle 1\rangle$
      and
      $\mathcal{I}_{\overline{a},x_1}(H)\ne\mathcal{I}_{\overline{a}}(H)$.
      In this case $H$ is polynomial but not strongly polynomial.
\end{enumerate}
\end{thm}
\begin{proof}
  First, assume $\mathcal{I}(H)=\langle 1\rangle$. In this case,
  $\mathcal{I}_{\overline{a}}(H)=\langle 1\rangle$, too. By the
  Nullstellensatz, $U(H)=\emptyset$. Then, for all standard polynomial
  $\Phi(x,y)$, we have $U(H,\Phi)=\emptyset$. Applying
  Theorem~\ref{HU}, we conclude that $G(\Phi)^*$ has no component
  isomorphic to $H$. Therefore, $H$ is not a polynomial graph.

  Now, assume $\mathcal{I}(H)\ne \langle 1\rangle$. By the
  Nullstellensatz we have
  $\mathbb{V}(\mathcal{I}_{\overline{a}}(H))\ne\emptyset$.  Note that
  no graph $H$ is $\Phi$-polynomial for all $\Phi(x,y)$.
  Then, by Proposition~\ref{proper}, there
  exists $\Phi(x,y)$ such that $U(H,\Phi)= \emptyset$. This implies
  $\mathcal{I}_{\overline{a}}\ne\{0\}$.

  It is convenient to label vertices
  $1,\ldots,n$ of $H$ in such a way that each vertex $i\ge
  2$ is adjacent to some vertex $j<i$. This can be done, for example
  by putting the labels on the vertices following the generation of a
  spanning tree by the Depth First Search (DFS) algorithm~\cite{KoVy}.

  In the second case,
  $\mathcal{I}_{\overline{a},x_1}(H)=\mathcal{I}_{\overline{a}}(H)$.
  Let $\overline{c}\in\mathbb{V}(\mathcal{I}_{\overline{a}}(H))$. If
  $\Phi_{\overline{c}}$ is non standard, Proposition~\ref{HUns}
  implies that $H=K_{n}$, a contradiction. Hence,
  $\Phi_{\overline{c}}(x,y)$ is standard.

  Let $u_1$ be a vertex of $G(\Phi_{\overline{c}})^*$.  Write
  $\Phi_{\overline{c}}(x,y)$ in the form
  $\Phi_{\overline{a}}(x,y)=\sum_{i=0}^d a_i(x)y^i$.  The hypothesis
  $\mathcal{I}_{\overline{a},x_1}(H)=\mathcal{I}_{\overline{a}}(H)$
  implies
  $(\overline{c},u_1)\in\mathbb{V}(\mathcal{I}_{\overline{a},x_1}(H))$.
  By induction, suppose that we have a partial solution
  $(\overline{c}, u_1, \ldots,
  u_k)\in\mathbb{V}(\mathcal{I}_{\overline{a}, x_1,\ldots, x_k}(H))$.
  Because of the labelling of the vertices of $H$, for some $j<k+1$,
  the vertex $u_j$ is adjacent to the vertex $u_{k+1}$. Then, the
  polynomial $\Phi(x_j,x_{k+1})$ belongs to
  $I_{\overline{a},x_1,\ldots,x_{k+1}}(H)$. Moreover, $a_d(u_j)\ne 0$
  because $u_j$ is a vertex of the non singular component
  $G(\Phi_{\overline{c}},u_1)$.  By the Extension
  Theorem~\cite{CoLiSh}, the partial solution extends to a solution
  $(\overline{c},u_1,\ldots,u_{k+1})$. Therefore, $(u_1,\ldots,u_n)$
  is a point of $U(H,\Phi_{\overline{c}})$. By
  Proposition~\ref{proper}, $\langle u_1,\ldots,u_n\rangle$ is a
  component of $G(\Phi_{\overline{c}})^*$ isomorphic to $H$.
  Therefore, $H$ is strongly $\Phi_{\overline{c}}$-polynomial for all
  $\overline{c}\in\mathbb{V}(\mathcal{I}_{\overline{a}}(H))$.

Finally, assume $\mathcal{I}_{\overline{a},x_1}(H)\ne
\mathcal{I}_{\overline{a}}(H)$.  As before, $U(H)\ne \emptyset$,
and, for any proper point $(\overline{c},u_1,\ldots,u_n)\in U(H)$
the graph $\langle u_1,\ldots,u_n\rangle$ is a component of
$G(\Phi_{\overline{c}})^*$. Therefore $H$ is
$\Phi_{\overline{c}}$-polynomial.  But as the indeterminate $x_1$
is not free, $H$ is not strongly $\Phi_{\overline{c}}$-polynomial.
\end{proof}

If $H=K_n$, then $\mathcal{I}(K_n)\ne \langle 1\rangle$ and
$\mathcal{I}_{\overline{a},x_1}(K_n)=\mathcal{I}_{\overline{a}}(K_n)$,
as in (ii). In this case, besides the standard polynomials
$\Phi_{\overline{c}}(x,y)$ such that $K_n$ is strongly
$\Phi_{\overline{c}}$-polynomial, there are points in
$\mathbb{V}(\mathcal{I}_{\overline{a}})$ corresponding to non
standard polynomials as described in Proposition~\ref{HUns}.

The proof of case (ii) provides also some insight about the
singular components:

\begin{prop}
Let $H$ be a strongly $\Phi$-polynomial graph. If $C$
is a singular component of $G(\Phi)$ without defective vertices, then
there exists an improper point $(u_1,\ldots,u_n)$ in $U(H,\Phi)$ such that
$C=\langle u_1,\ldots,u_n\rangle$. Moreover there exists a graph morphism
from $H$ onto $C$.
\end{prop}
\begin{proof}
  In order to apply the Extension Theorem to
  Theorem~\ref{characideal} (ii), the crucial point is the condition
  $a_d(u_j)\ne 0$, that means that $u_j$ is not a defective vertex.
  Therefore if $u_1$ is taken in a singular component $C$ without
  defective vertices, then $(u_1,\ldots,u_n)$ is an improper point in
  $U(H,\Phi)$. It is easily checked that $i\mapsto u_i$ is an
  exhaustive morphism from $H$ to $C$.
\end{proof}

Conjecture~\ref{ccc} is equivalent to saying that  case (iii) in
Theorem~\ref{characideal} never occurs. The following proposition
reduces the scope of the conjecture.

\begin{prop}
\label{rsccc}
Let $\Phi(x,y)$ be a standard symmetric polynomial of partial degree
$d$. Then,
\begin{enumerate}
\item[\emph{(i)}] If $G(\Phi)^*$ has an uncountable number of finite components, then all components
of $G(\Phi)^*$ are finite and isomorphic.
\item[\emph{(ii)}] If $G(\Phi)^*$ has a countable number of finite components, then each finite component
of $G(\Phi)^*$ is not a strongly polynomial graph.
\end{enumerate}
\end{prop}
\begin{proof}
(i) There exist a countable number of non isomorphic finite
graphs. As
  the number of components of $G(\Phi)^*$ is uncountable,
  a finite graph $H$ exists such that an uncountable number of components of
  $G(\Phi)^*$ are isomorphic to $H$. Then $U(H,\Phi)$ is 1-dimensional
  and, for each vertex $u_1$ of $G(\Phi)^*$, we have a proper point
  of $U(H,\Phi)$ with $x_1=u_1$, i.e. a component isomorphic to $H$.

  (ii) Let $H$ be a finite component of $G(\Phi)^*$. Consider the
  three cases of Theorem~\ref{characideal}. As $H$ is a
  $\Phi$-polynomial graph, we are not in case (i). As there are
  only a countable number of finite components, we are not in case
  (ii). Then, case (iii) applies.
\end{proof}

Thus, Conjecture~\ref{ccc} is reduced to the following: there
exists no standard symmetric polynomial $\Phi(x,y)$ such that
$G(\Phi)^*$ has a countable number of finite components any of
which is isomorphic to a strongly polynomial graph.

On the other side, computational evidence suggests that if
$\Phi(x,y)$ is a standard symmetric polynomial and $G(\Phi)^*$ has
infinite graphs as components, then it is not true that all
components of $G(\Phi)^*$ are isomorphic. For instance, this seems
to be the case with the polynomial $\Phi(x,y)=x^3+y^3+xy-1$.

\section{General algorithm}
\label{GenAlg} Given a finite, connected, $d$-regular graph $H$,
we want to determine its characteristic ideal
$\mathcal{I}_{\overline{a}}(H)=\mathbb{I}(U(H))\cap
\mathbb{C}[\overline{a}]$. We start with the system of polynomials
$S(H)$, choose the monomial order $
\operatorname{lex}(x_n,\dots,x_1,a_{0\,0},a_{1\,0},a_{1\,1},\dots,a_{d\,d}),
$ and use the generalized gaussian elimination algorithm {\tt gge}
in the \emph{Maple} library {\tt dpgb}~\cite{Montes} in order to
simplify $S(H)$. At any step before launching Buchberger's
algorithm, we must eliminate factors of the form $x_i-x_j$ in
every new polynomial generated. Reductions and Buchberger's
algorithm can be combined, to obtain the Gr\"obner basis of the
ideal $\mathcal{I}(H)$. The polynomials in this basis depending
only on the variables $\overline{a}$, are the Gr\"obner basis of
the characteristic ideal $\mathcal{I}_{\overline{a}}(H)$. Then we
can also test if
$\mathcal{I}_{\overline{a},x_1}(H)=\mathcal{I}_{\overline{a}}(H)$
to decide, by Theorem~\ref{characideal}, if $H$ is strongly
polynomial.

To make the computation effective it is strictly needed to add to
$S(H)$ as many polynomials in $\mathbb{I}(U(H))\setminus
\mathbb{I}(Z(H))$ as possible.  Before giving a method for
obtaining polynomials of this kind, let us consider an example.
Let $H$ be the 4-cycle $C_4$, and consider the system
$$
S(C_4)=\{\Phi(x_1,x_2),\Phi(x_2,x_3,\Phi(x_3,x_4),\Phi(x_4,x_1)\},
$$
where $\Phi=\Phi_{\overline{a}}$. For a given
$\overline{a}\in\mathbb{C}^m$ and $u \in \mathbb{C}$ let
$\lambda_1,\lambda_2$ be the two roots of $\Phi(u,y)$.  Then
$(\overline{a}, u,\lambda_1,u,\lambda_2)$,
$(\overline{a},u,\lambda_2,u,\lambda_1)$,
$(\overline{a},u,\lambda_1,u,\lambda_1)$ and
$(\overline{a},u,\lambda_2,u,\lambda_2)$ are improper points in
$W(C_4)$.  Let now $\mu_1,\mu_2$ be such that $u,\mu_i$ are the
two distinct roots of $\Phi(\lambda_i,y)$, for $i=1,2$. If
$\mu_1=\mu_2$, then $(\overline{a},u,\lambda_1,\mu_1,\lambda_2)$
and $(\overline{a},u,\lambda_2,\mu_1,\lambda_1)$ are proper points
in $U(C_4)$.  Thus, for any $\overline{a}$ and $u$ there exist a
finite number of solutions in $Z(C_4,\Phi_{\overline{a}})$, and
this variety is of dimension 1 for any $\overline{a}$. But the
condition $\mu_1=\mu_2$ will be satisfied only if $C_4$ is a
$\Phi_{\overline{a}}$-polynomial graph, and there are proper
points in $U(C_4,\Phi_{\overline{a}})$. If $\mu_1=\mu_2$ for any
$u$, then $C_4$ is strongly $\Phi_{\overline{a}}$ -polynomial. The
undesired solutions in $Z(C_4)$ appear owing to the fact that  no
distinction is made in $S(C_4)$ between the $y$-roots of
$\Phi(u,y)$.

Let $H$ be a finite, connected, $d$-regular graph.  The following
method allows to obtain a set of polynomials in $\mathbb{I}(U(H))$
(depending on a vertex of $H$) that separates roots. Consider a
vertex $i_0$ in $H$ and let $x_{i_1},\ldots x_{i_d}$ be the
indeterminates corresponding to the vertices adjacent to $i_0$. In
the following discussion we write $x_j$ instead of $x_{i_j}$ to
avoiding subscripts. Consider the polynomials $\Phi(x_0,x_j)$,
$1\le j\le d$, in $S(H)$.   Set
$$
\Phi_0(x_0;x_1)=\Phi(x_0,x_1),
$$
and define recursively
\[
\Phi_{\ell-1}(x_0;x_1,\dots,x_{\ell})=
  \displaystyle\frac{\Phi_{\ell-2}(x_0;x_1\dots,x_{\ell-1})
        -\Phi_{\ell-2}(x_0;x_1\dots,x_{\ell-2},x_{\ell})}{x_{\ell}-x_{\ell-1}}.
\]

\begin{prop}
\label{improve}
The polynomials $\Phi_{\ell}$ have the following properties:
\begin{enumerate}
  \item [\emph{(i)}] $\Phi_{\ell-1}$ is a polynomial of $U(H)$ for $1\le \ell \le d$.
  \item [\emph{(ii)}] $\Phi_{\ell-1}(x_0;x_1,\dots,x_{\ell})$ is symmetrical
    in the set of variables $\{x_1,\dots,x_{\ell}\}$, and its total
    degree as a polynomial in the variables $\overline{x}$ is $2d-\ell+1$.
  \item [\emph{(iii)}] If $\Phi(u_0,y)$ has $d$ different roots
  $u_1,\dots,u_{d}$, then the set of solutions of the
  system formed by the $d$ polynomials
$\Phi_0(u_0;x_1),\  \Phi_1(u_0;x_1,x_2), \ \dots,\
\Phi_{d-1}(u_0;x_1,\dots,x_{d})$ is exactly the set of all
permutations of the solution $\{x_1=u_1, \ldots, x_d=u_{d}\}$.
\end{enumerate}
\end{prop}
\begin{proof}
(i) For $\ell=2$ we have
\[
\Phi_0(x_0;x_1)-\Phi_0(x_0;x_2)=
 \displaystyle\sum_{i,j=0}^d a_{ij} x_0^i (x_1^j-x_2^j)
 \ \ =\displaystyle\sum_{i,j=0
 }^d a_{ij} x_0^i (x_1-x_2)\sum_{k=1}^j x_1^{j-k}
x_2^{j-1}.
\]
Thus $\Phi_1(x_0;x_1,x_2)$ belongs to $U(H)$ and we have
\[
\Phi_1(x_0;x_1,x_2)=
 \displaystyle\sum_{i,j=0}^d a_{ij} x_0^i \sum_{k=1}^j x_1^{j-k}
x_2^{j-1}.
\]
Iterating, we obtain an explicit formula for $\Phi_{\ell-1}$:
\[
\label{recurexpl} \Phi_{\ell-1}(x_0;x_1,\dots,x_{\ell})=
   \displaystyle\sum_{i,j=0}^d a_{ij} x_0^i \sum_{k_1=1}^j x_1^{j-k_1}
\sum_{k_2=1}^{k_1}x_2^{k_1-k_2} \dots
   \displaystyle\sum_{k_{\ell-1}=1}^{k_{\ell-2}}
x_{\ell-1}^{k_{\ell-2}-k_{\ell-1}} x_{\ell}^{k_{\ell-1}-1}.
\]
showing that it belongs to $U(H)$.

(ii) It can be proved by induction that formula~(\ref{recurexpl})
is equivalent to
$$
\begin{array}{l}
\Phi_{\ell-1}(x_0;x_1,\dots,x_{\ell})=
  \displaystyle\sum_{i,j=0}^d a_{ij} x_0^i \sum_{\overline{k}}
 x_1^{k_1}\cdots
x_{\ell}^{k_{\ell}},  \\
\end{array}
$$
where the sum over $\overline{k}$ is extended to all
$\overline{k}=(k_1,\dots,k_{\ell})$ verifying $k_i \ge 0$
and $\sum_{s=1}^{\ell}k_s =j-\ell+1$.  This formula is
explicitly symmetric in the set of variables $\{x_1,\dots,x_{\ell}\}$,
and its degree in $\overline{x}$ is obviously $2d-\ell+1$.

(iii) $\Phi(u_0,x_1)$ has exactly the $d$ solutions
  $\{x_1=u_1,\dots,x_1=u_d \}$. Then $\Phi_0(u_1;u_i)=0$ and
  $\Phi_1(u_0;u_i,x_2)=0$ imply $\Phi_0(u_i;x_2)=0$ and thus
  $\Phi_1(u_0;u_i,x_2)$ has the same roots as $\Phi(u_0,x_2)$ except
  for $u_i$. Similarly, we can prove that
  $\Phi_{\ell-1}(u_0;u_{i_1},\dots,u_{i_{\ell-1}},x_{\ell})$ has the same roots
  as $\Phi(u_0,x_\ell)$ except for $\{ u_{i_1},\dots,u_{i_{\ell-1}}\}$.  Thus
  the set of solutions of $S(u_0)$ is the set of all
  permutations of $\{x_1=u_1,\dots,x_d=u_d\}.$
\end{proof}

Let $V_i$ be the set of vertices of $H$ adjacent to the vertex $i$.
The \emph{completed system} $S'(H)$ is formed by all the polynomials
$\Phi_{\ell-1}(x_0;x_{i_1},\ldots,x_{i_\ell})$, where
$\{i_1,\ldots,i_{\ell}\}$ is a $\ell$-subset of $V_i$, for all
$\ell\in[d]$ and $i\in[n]$. Note that for $\ell=1$ we obtain the
polynomials in $S(H)$.
The
number of polynomials in $S'(H)$ is
$$
n\sum_{\ell=1}^d {d \choose \ell}= n(2^d-1).
$$
Nevertheless, in this account there are repeated polynomials. For
instance $\Phi(x_i,x_j)=\Phi_0(x_i;x_j)$ appear twice. The system
$S'(H)$ being a set, repetitions have to be crossed out.

In general, the solutions of the completed system  $S'(H)$ are not
exactly the points in $U(H)$. Factors $x_i-x_j$ can appear in the
computing of a Gr\"obner basis. Often, it is possible to take into
account the symmetry of the graph in order to eliminate these
extraneous solutions, by introducing a new set of reduced
polynomials. For instance, for graphs with cliques of order
$\ell+1$ the following polynomials are helpful. Set
$$
\Phi_{\ell-1,0}(x_0;x_1,\dots,x_{\ell})
=\Phi_{\ell-1}(x_0;x_1,\dots,x_{\ell}),
$$
and, recursively,
\begin{eqnarray*}
\Phi_{\ell-1,k}(x_0,\dots,x_k;x_{k+1},\dots,x_{\ell})&=&
   \left(\Phi_{\ell-1,k-1}(x_0,\dots,x_{k-2},x_k;
                   x_{k-1},x_{k+1}\dots,x_{\ell})\right.\\
    && \left. {} - \Phi_{\ell-1,k-1}(x_0,\dots,x_{k-1};x_{k},\dots,x_{\ell})\right) /(x_{k-1}-x_k).
\end{eqnarray*}

\begin{prop}
\label{polycliques}
Suppose that the vertices $\{0,1,\dots,\ell\}$ is a clique of $H$. Then,
\begin{itemize}
  \item [\emph{(i)}] $\Phi_{\ell-1,k}$ are polynomials in $\mathbb{I}(U(H))$ for $0\le k \le
  \ell-1$.
  \item [\emph{(ii)}] $\Phi_{\ell-1,k}$ are symmetric in the second set of
  variables.
  \item [\emph{(iii)}] The degree of $\Phi_{\ell-1,k}$ in $\overline{x}$ is $2d-\ell+1-k$.
\end{itemize}
\end{prop}
\begin{proof}
The fact that each two vertices in $0,\ldots,\ell$ are adjacent
implies that $\Phi_{\ell-1,k}$ are polynomials of
$\mathcal{I}(H)$. The symmetry of the
 $\Phi_{\ell-1}$ in the second set of variables produces the symmetry of the
 $\Phi_{\ell-1,k}$ in the second set of variables and their
 degree in $\overline{x}$ is deduced directly from the degree of the $\Phi_{\ell-1}$.
\end{proof}

\subsection{Application}
We apply the general algorithm to determine the characteristic
ideal of a 3-cycle and of a 4-cycle. Let
\begin{eqnarray*}
\Phi(x,y)&=&a_{00}+a_{10} (x+y)+a_{11} x y+a_{20} (x^2+y^2)
          +a_{21}xy(x+y)+a_{22} x^2 y^2,\\
\Phi_1(x;y,z)&=&  a_{10}+a_{11} x+a_{20} (y+z)
              {}+a_{21} x (x+y+z)+a_{22} x^2 (y+z),\\
\Phi_{11}(x,y;z)&=&a_{11}-a_{20}+a_{21} (x+y+z)
                 +a_{22}(xy+yz+zx).
\end{eqnarray*}
The complete system for the
3-cycle is
\[
 S'(C_3) =\{\Phi(x_1,x_2),\Phi(x_2,x_3),\Phi(x_3,x_1),
 \Phi_1(x_1;x_2,x_3),\Phi_1(x_2;x_3,x_1),\Phi_1(x_3;x_1,x_2)\}.\\
\]
As $C_3$  is a complete graph, we add the polynomial
$\Phi_{11}(x_1,x_2;x_3)$. Let
$S''(H)=S'(H)\cup\{\Phi_{11}(x_1,x_2;x_3)\}$. We take the monomial
order
$$
\operatorname{lex}(x_3,x_2,x_1,a_{00},a_{10},a_{11},a_{20},a_{21},a_{22})
$$
and calculate the Gr\"obner basis of $S''(C_3)$, which is easily
computed and contains 9 polynomials. The quick computation is
owing to the inclusion of the polynomial $\Phi_{11}$, that
reflects the symmetry.  The Gr\"obner basis provides the following
elimination ideal:
$$
\mathcal{I}_{\overline{a}}(C_3)=
\langle {\mathbf{a}_{00}\,a_{22}}+a_{20}\,a_{11}-a_{20}^2-a_{21}\,a_{10}\rangle
$$
and $\mathcal{I}_{\overline{a}\,x_1}(C_3)=\mathcal{I}_{\overline{a}}(C_3)$.

Consider now $4$-cycles. The complete system is:
\begin{eqnarray*}
 S'(C_4)&=&\{\Phi(x_1,x_2),\Phi(x_2,x_3),\Phi(x_3,x_4),\Phi(x_4,x_1),\\
 &&\phantom{\{} \Phi_1(x_1;x_4,x_2),\Phi_1(x_2;x_1,x_3),\Phi_1(x_3;x_2,x_4),\Phi_1(x_4;x_3,x_1)\}.\\
\end{eqnarray*}
The computations become only effective when we add a
new reduced polynomial that reflects the symmetry, and eliminates
the extraneous solution $x_1=x_3$, namely:
\begin{eqnarray*}
\lefteqn{\Psi(x_1,x_3;x_2,x_4)
=\frac{\Phi_1(x_1;x_2,x_4)-\Phi_1(x_3;x_2,x_4)}{x_1-x_3}}\\
&=& a_{11}+a_{21}(x_1,x_2,x_3,x_4)+a_{22}(x_1+x_3)(x_2+x_4).
\end{eqnarray*}
Take $S''(C_4)=S'(C_4)\cup \{\Psi(x_1,x_3;x_2,x_4)\}$. The direct
computation of the Gr\"obner basis, when using an automatic
method, becomes difficult. We use the technique of stopping the
computation when a high number of polynomials have been computed
and then use {\tt gge} routine in the {\tt dpgb} library to reduce
the basis. The result is a basis of 24 polynomials, which provides
the following characteristic ideal:
\[
\mathcal{I}_{\overline{a}}(C_4)=
  \langle{\bf a_{00} a_{11} a_{22}}+2 a_{21} a_{20} a_{10}
 -a_{20}^2 a_{11} -a_{21}^2 a_{00}-a_{10}^2 a_{22}\rangle.
\]
and, as before,
$\mathcal{I}_{\overline{a} x_1}(C_4)=\mathcal{I}_{\overline{a}}(C_4)$.

As shown in~\cite{BrMo}, the polynomial of partial degree two
$\Phi_{\overline{a}}(x,y)$ can be reduced by a translation to a
polynomial with $a_{21}=0$. By performing the above computations
in this case, the number of polynomials in the basis reduces to 8
polynomials for $\mathbb{I}(S''(C_3))$ and 19 for
$\mathbb{I}(S''(C_4))$.

\section{Cycles}
\label{AlgCompCycles}

In~\cite{BrMo} a complete study of the components of $G(\Phi)$
when $\Phi(x,y)$ is a symmetric polynomial of total degree two is
given. The method can be used to determine conditions on the
coefficients of a polynomial $\Phi(x,y)=a(x)y^2+b(x)y+c(x)$ of
partial degree 2 for obtaining cycles of length $n$ as components
of $G(\Phi)^*$. Let
\begin{eqnarray*}
a(x) &=& a_{22}x^2+a_{21}x+a_{20}, \\
b(x) &=& a_{21}x^2+a_{11}x+a_{10},\\
c(x) &=& a_{20}x^2+a_{10}x+a_{00}.
\end{eqnarray*}
As a polynomial in $y$, the  sum of the two roots of $\Phi(x,y)$
equals  $-b(x)/a(x)$. Then, we have the recurrence:
$$
v_n=-v_{n-2}-\frac{b(v_{n-1})}{a(v_{n-1})}=\frac{p_n}{q_n}.
$$
By iterating the recurrence with free initial values $v_0$ and
$v_1$, we obtain, by substitution and simplification, expressions
for $p_n$ and $q_n$, in terms of $v_0$, $v_1$ and of the
coefficients $\overline{a}$. To obtain $n$-cycles we must impose
$K_n=p_n -v_0 q_n=0$ and $\Phi(v_0,v_1)=0$. We use the above
conditions, dividing $K_n$ by $[\Phi(v_0,v_1)]$ using a convenient
monomial order. The result is a polynomial that has one factor
depending only on the parameters $\overline{a}$. Consequently the
polynomial produces $n$-cycles for any initial point $v_0$, when
the factor containing only the parameters vanishes. In this way we
obtain the characteristic ideals for 3, 4 and 5 cycles, which are
principal ideals. These are:
\begin{eqnarray*}
\Delta_3 &=& {\bf a_{22}\, a_{00}}+a_{11}\, a_{20}-a_{20}^2\,-a_{21} a_{10},\\
\Delta_4 &=& {\bf a_{22}\, a_{11}\, a_{00}}-a_{22}\,
                  a_{10}^2\,-a_{11}\, a_{20}^2+2 a_{21}\, a_{20}\, a_{10}
              -a_{21}^2\, a_{00},\\
\Delta_5 &=& {\bf  a_{22}^3\, a_{00}^3}-a_{21}^3\, a_{10}^3-4 a_{22}\,
                   a_{20}^3\, a_{10}^2+5\,a_{21}\, a_{20}^4\, a_{10}
             +a_{20}^2 \,a _{10}^2\, a_{21}^2 \\
         &&         +a_{10}^2\, a_{21}^2\, a_{20}\, a_{11}
                  -4\, a_{10}\, a_{21}\, a_{20}^3\, a_{11}
               -a_{22}^2\, a_{10}^4-a_{11}\, a_{20}^5 -a_{20}^6 \\
          &&        +3\, a_{22}\, a_{11}\, a_{10}^2\, a_{20}^2
             {} +a_{22}\, a_{21}\,a_{11}\, a_{10}^3
                  -a_{22}\, a_{11}^2\, a_{20}\, a_{10}^2+a_{11}^2 \,a_{20}^4\\
         &&    {} +4\, a_{22}^2\, a_{20}\, a_{10}^2\, a_{00}
                  +3\, a_{11}\, a_{21}^2\, a_{20}^2\, a_{00}
              -2\, a_{22}\, a_{21}\, a_{20}^2\, a_{10}\, a_{00}\\
          &&        -a_{20}\, a_{11}^2\, a_{21}^2\, a_{00}
             {}-4\, a_{20}^3\, a_{21}^2\, a_{00}
                 -3\, a_{22}^2\, a_{20}^2\, a_{00}^2
             {}+a_{22}^2\, a_{10}^2\, a_{11}\, a_{00}\\
           &&       +a_{11}\, a_{21}^2\, a_{22}\, a_{00}^2
             {} -3\, a_{22}^2\, a_{10}\, a_{21}\, a_{00}^2
             {}+a_{10}\, a_{11}\, a_{21}^3\, a_{00}
                 -a_{22}\, a_{20}^2\, a_{11}^2\, a_{00}\\
         &&    {}-4\, a_{22}\, a_{10}\, a_{11}\, a_{20}\, a_{21}\, a_{00}
                 +a_{22}\, a_{20}\, a_{11}^3\, a_{00}
             {}+2\, a_{22}\, a_{11}\, a_{20}^3\, a_{00}\\
          &&        +4\, a_{22}\, a_{21}^2\, a_{20}\, a_{00}^2
             {}-a_{11}\, a_{22}^2\, a_{20}\, a_{00}^2
                  +a_{22}\, a_{10}^2\, a_{21}^2\, a_{00}\\
         &&    {}-a_{22}\, a_{10}\, a_{21}\, a_{11}^2\, a_{00}
                 -a_{21}^4\, a_{00}^2 +3 a_{22}\, a_{20}^4\, a_{00}.
\end{eqnarray*}
Using the above characteristics ideals, it is easy to obtain
examples of polynomials producing cycles:\\

\centerline{
\begin{tabular}{r|l}
Graph & Polynomial \\
\hline
 $C_3$ & $ x^2 y^2+x^2+y^2-x y+2$ \\
 $C_4$ & $ x^2 y^2+x^2+y^2+x y+1$\\
 $C_5$ & $ x^2 y^2+x^2+y^2-2 x y+x+y-2$.\\
\end{tabular}
}

\section{Complete graphs}
\label{AlgComplGraph} For complete graphs $K_{d+1}$ we use a
specific technique that takes into account the symmetry of the
graph. We start writing the system $S(H)$ of polynomials
corresponding to $K_{d+1}$. Then, as the number of parameters
$\overline{a}$ is $(d+2)(d+1)/2$, and the number of edges (=
equations) is $d(d-1)/2$ we can solve the linear system
considering the $\overline{a}$ as variables. This provides some of
the $\overline{a}$ in terms of the rest. In order to obtain the
correct result, it is important to choose the coefficients with
greatest indices as parameters and to express the $\overline{a}$
with smaller indexes in terms of them. Being careful we can obtain
an expression for some of the $\overline{a}$ linearly dependent in
the rest of the $\overline{a}$, and polynomial in the
$\overline{x}$. Owing to the symmetry of the complete graph in the
vertices, we can now transform the dependence of these expressions
in the $\overline{x}$ in terms of the elementary symmetrical
polynomials of the $\overline{x}$ say $s_1,s_2,\dots,s_d$.

The resulting system of equations turns out to be linear in the
$s_i$ and very simple. For $K_3$, $K_4$ and $K_5$ the
corresponding set of polynomials defining the systems are:
$$
\begin{array}{l}
S(K_3) = \left\{
\begin{array}{l}
 a_{00}-a_{2 0}\, s_{2}+a_{2 1}\,(-1+s_{1}), \\
 a_{1 0}+a_{2 0}\,s_{1}+a_{2 2}\,(-1+s_{1}), \\
 a_{1 1}-a_{2 0}+a_{2 1}\,s_{1}+a_{2 2}\,s_{2}.\\
\end{array}
\right.
\\
S(K_4) = \left\{
\begin{array}{l}
 a_{00}+a_{3 0}\,s_{3}+a_{3 1}\,s_{1}, \\
 a_{1 0}-a_{3 0}\,s_{2}+a_{3 2}\,s_{1} ,\\
 a_{1 1}+a_{3 0}\,s_{1}-a_{3 1}\,s_{2}-a_{3 2}\,s_{3}+a_{3 3}\,s_{1}, \\
 a_{2 1}-a_{3 0}+a_{3 1}\,s_{1}-a_{3 3}\,s_{3}, \\
 a_{2 2}-a_{3 1}+a_{3 2}\,s_{1}+a_{3 3}\,s_{2}, \\
 a_{2 0}+a_{3 0}\,s_{1}+a_{3 3}\,s_{1}. \\
\end{array} \right.
\\
S(K_5) = \left\{
\begin{array}{l}
  a_{00}-a_{4 0}\,s_{4}+a_{4 1}\,(-1+s_{1}), \\
  a_{1 0}+a_{4 0}\,s_{3}+a_{4 2}\,(-1+s_{1}), \\
  a_{1 1}-a_{4 0}\,s_{2}+a_{4 1}\,s_{3} +a_{4 2}\,s_{4}+a_{4 3}\,(-1+s_{1}), \\
  a_{2 0}-a_{4 0}\,s_{2}+a_{4 3}\,(-1+s_{1}), \\
  a_{2 1}+a_{4 4}\,(-1+s_{1})+a_{4 0}\,s_{1}-a_{4 1}\,s_{2}+a_{4 3}\,s_{4}, \\
  a_{2 2}+a_{4 4}\,s_{4}-a_{4 0}+a_{4 1}\,s_{1}-a_{4 2}\,s_{2}-a_{4 3}\,s_{3}, \\
  a_{3 0}+a_{4 4}\,(-1+s_{1})+a_{4 0}\,s_{1}, \\
  a_{3 1}+a_{4 4}\,s_{4}-a_{4 0}+a_{4 1}\,s_{1}, \\
  a_{3 2}-a_{4 4}\,s_{3}-a_{4 1}+a_{4 2}\,s_{1}, \\
  a_{3 3}+a_{4 4}\,s_{2}-a_{4 2}+a_{4 3}\,s_{1} . \\
\end{array}
\right.
\end{array}
$$

As we see, the equations do not depend on $s_d$, the latest
ele\-men\-ta\-ry symmetrical polynomial. This proves directly the
conjecture, namely the ideal $\mathcal{I}$ in the variables
$s_1,\dots,s_d,\overline{a}$ has one degree of freedom more than
the elimination ideal in the variables $\overline{a}$, and the
variable $s_d$ is free.

Now we apply the standard method with the new variables
$\overline{s}$, using the order $ \succ_s
=\operatorname{lex}(s_1,s_2,\dots,s_{d-1},a_{00},a_{10},\dots,a_{dd}),
$ and determine the Gr\"obner basis of the ideals
$\mathcal{I}_{\overline{a}}(K_n)$. In this way, we obtain the
characteristic ideals for $K_3$, $K_4$, $K_5$ and $K_6$. These are
$$
\begin{array}{lll}
\mathcal{I}_{\overline{a}}(K_3) &=& \langle {\bf a_{00}\,a_{22}}+a_{20}\,a_{11}-a_{20}^2-a_{21}\,a_{10}\rangle.\\
&& \\
\mathcal{I}_{\overline{a}}(K_4) &=&\langle{\bf a_{1 1}\,a_{3
3}}-a_{3 2}\,a_{2 1}+a_{3 2}\,a_{3 0}-a_{31}^2
      +a_{3 1}\,a_{2 2}-a_{3 3}\,a_{2 0}, \\
  && {\bf a_{1 0}\,a_{3 3}}-a_{3 2}\,a_{2 0}
   -a_{3 0}\,a_{3 1}+a_{30}\,a_{2 2}, \\
   && {\bf a_{1 0}\,a_{2 1}\,a_{3
2}}-a_{1 0}\,a_{3 1}\,a_{2 2}
 {} -a_{1 0}\,a_{3 2}\,a_{3 0} +a_{1 0}\,a_{3 1}^2
 -a_{1 1}\,a_{3 2}\,a_{2 0}\\
 && +a_{1 1}\,a_{3 0}\,a_{2 2} -a_{1
1}\,a_{3 0}\,a_{3 1} {} +a_{3 2}\,a_{2 0}^2
 -a_{2 0}\,a_{3 0}\,a_{2 2}+a_{2 0}\,a_{30}\,a_{3 1},\\
&&
{\bf a_{0 0}\,a_{3 3}}-a_{3 1}\,a_{2 0}+a_{3 0}\,a_{2 1}-a_{30}^2,\\
&&
{\bf a_{0 0}\,a_{3 2}}-a_{3 1}\,a_{1 0}+a_{3 0}\,a_{1 1}-a_{30}\,a_{2 0},\\
&& {\bf a_{0 0}\,a_{2 2}}-a_{1 0}\,a_{2 1}+a_{1 0}\,a_{3 0}
     -a_{20}^2+a_{2 0}\,a_{1 1}-a_{3 1}\,a_{0 0}\rangle.\\ &&
\\
\mathcal{I}_{\overline{a}}(K_5) &=& \langle  {\bf a_{2 2}\,a_{4
4}}-a_{4 4}\,a_{31}+a_{4 1}\,a_{4 3}
 -a_{4 3}\,a_{3 2}-a_{4 2}^2+a_{4 2}\,a_{3 3}, \\
&& {\bf a_{2 1}\,a_{4 4}}-a_{4 4}\,a_{3 0}+a_{4 3}\,a_{4 0}
 -a_{4 3}\,a_{3 1}-a_{4 1}\,a_{4 2}+a_{4 1}\,a_{3 3}, \\
&& {\bf a_{2 1}\,a_{3 2}\,a_{4 3}}+a_{4 3}\,a_{4 1}\,a_{3 0}-a_{3
0}\,a_{4 3}\,a_{3 2}
 -a_{3 0}\,a_{4 2}^2
 +a_{3 0}\,a_{4 2}\,a_{3 3}-a_{4 3}\,a_{4 0}\,a_{3 1}\\
 && +a_{4 3}\,a_{3 1}^2
 +a_{3 1}\,a_{4 1}\,a_{4 2}
 -a_{3 1}\,a_{4 1}\,a_{3 3}
 -a_{2 1}\,a_{4 1}\,a_{4 3}+a_{2 1}\,a_{4 2}^2-a_{2 1}\,a_{4 2}\,a_{3 3}\\
&&{} +a_{2 2}\,a_{4 3}\,a_{4 0}-a_{2 2}\,a_{4 3}\,a_{3 1}
 -a_{2 2}\,a_{4 1}\,a_{4 2}+a_{2 2}\,a_{4 1}\,a_{3 3}, \\
&&
{\bf a_{2 0}\,a_{4 4}}-a_{4 3}\,a_{3 0}-a_{4 2}\,a_{4 0}+a_{4 0}\,a_{3 3}, \\
\end{array}
$$
$$
\begin{array}{lll}
&& {\bf a_{2 0}\,a_{3 2}\,a_{4 3}}-a_{4 2}\,a_{3 3}\,a_{2 0}
 -a_{4 1}\,a_{4 3}\,a_{2 0}+a_{4 2}^2\,a_{2 0}
 -a_{2 2}\,a_{4 3}\,a_{3 0}+a_{2 2}\,a_{4 0}\,a_{3 3} \\
 && -a_{4 2}\,a_{4 0}\,a_{2 2}+a_{3 1}\,a_{4 3}\,a_{3 0}
  -a_{3 1}\,a_{4 0}\,a_{3 3}+a_{4 2}\,a_{4 0}\,a_{3 1},\\
&& {\bf a_{2 0}\,a_{3 1}\,a_{4 3}}-a_{4 1}\,a_{3 3}\,a_{2 0}+a_{4
1}\,a_{4 2}\,a_{2 0}
 -a_{2 1}\,a_{4 3}\,a_{3 0}
 +a_{2 1}\,a_{4 0}\,a_{3 3}+a_{4 3}\,a_{3 0}^2\\
 && -a_{3 0}\,a_{4 0}\,a_{3 3}
 +a_{4 0}\,a_{4 2}\,a_{3 0}
 -a_{4 0}\,a_{4 2}\,a_{2 1}-a_{4 0}\,a_{4 3}\,a_{2 0}, \\
&& {\bf a_{2 0}\,a_{3 1}\,a_{4 2}}+a_{2 2}\,a_{4 1}\,a_{3 0}-a_{4
2}\,a_{3 0}\,a_{2 1}
 -a_{3 2}\,a_{4 1}\,a_{2 0}
 -a_{3 1}\,a_{4 1}\,a_{3 0}+a_{4 1}^2\,a_{2 0}\\
 && +a_{4 2}\,a_{3 0}^2-a_{4 0}\,a_{3 2}\,a_{3 0}
 +a_{4 0}\,a_{3 2}\,a_{2 1}+a_{4 0}\,a_{3 1}^2-a_{2 2}\,a_{4 0}\,a_{3 1}
 -a_{4 0}\,a_{4 2}\,a_{2 0}\\
&&{}
 -a_{4 0}^2\,a_{3 1}+a_{4 0}^2\,a_{2 2}+a_{4 0}\,a_{4 1}\,a_{30}-a_{4 0}\,a_{4 1}\,a_{2 1}, \\
&& {\bf a_{1 1}\,a_{4 4}}-a_{4 3}\,a_{3 0}-a_{4 2}\,a_{3 1}
 -a_{4 1}^2+a_{4 1}\,a_{3 2}+a_{4 0}\,a_{3 3}, \\
&& {\bf a_{1 1}\,a_{4 3}}-a_{4 1}\,a_{3 1}+a_{4 1}\,a_{2 2}+a_{4
2}\,a_{3 0}
 -a_{4 2}\,a_{2 1}-a_{4 3}\,a_{2 0}, \\
&& {\bf a_{1 1}\,a_{3 3}}+a_{3 2}\,a_{3 0}-a_{3 2}\,a_{2 1}-a_{3
1}^2
 +a_{3 1}\,a_{2 2}-a_{4 2}\,a_{1 1}
 -a_{3 3}\,a_{2 0}+a_{4 2}\,a_{2 0}\\
 && +a_{4 0}\,a_{3 1}
 -a_{4 0}\,a_{2 2}-a_{4 1}\,a_{3 0}+a_{4 1}\,a_{2 1}, \\
&& {\bf a_{1 0}\,a_{4 4}}-a_{4 2}\,a_{3 0}+a_{4 0}\,a_{3 2}
 -a_{4 1}\,a_{4 0}, \\
&& {\bf a_{1 0}\,a_{4 3}}-a_{4 2}\,a_{2 0}+a_{4 0}\,a_{2 2}
 -a_{4 0}\,a_{3 1}, \\
&& {\bf a_{1 0}\,a_{3 3}}+a_{4 1}\,a_{2 0}-a_{1 0}\,a_{4 2}
 -a_{3 2}\,a_{2 0}-a_{3 0}\,a_{3 1}+a_{3 0}\,a_{2 2}, \\
&& {\bf a_{1 0}\,a_{3 1}\,a_{4 2}}-a_{3 2}\,a_{4 1}\,a_{1 0}
  -a_{4 0}\,a_{1 0}\,a_{4 2}+a_{4 1}^2\,a_{1 0}
 -a_{3 0}\,a_{4 2}\,a_{1 1}+a_{3 2}\,a_{4 0}\,a_{1 1} \\
 && -a_{4 1}\,a_{4 0}\,a_{1 1}+a_{3 0}\,a_{4 2}\,a_{2 0}
 -a_{4 0}\,a_{3 2}\,a_{2 0}+a_{4 0}\,a_{4 1}\,a_{2 0}, \\
&& {\bf a_{1 0}\,a_{2 1}\,a_{4 2}}-a_{2 2}\,a_{4 1}\,a_{1 0}
 -a_{3 0}\,a_{1 0}\,a_{4 2}+a_{3 1}\,a_{4 1}\,a_{1 0}
 -a_{2 0}\,a_{4 2}\,a_{1 1} +a_{2 2}\,a_{4 0}\,a_{1 1}\\
 &&-a_{3 1}\,a_{4 0}\,a_{1 1}+a_{4 2}\,a_{2 0}^2
 -a_{2 0}\,a_{4 0}\,a_{2 2}+a_{2 0}\,a_{4 0}\,a_{3 1}, \\
&& {\bf a_{1 0}\,a_{2 1}\,a_{3 2}}-a_{1 0}\,a_{4 0}\,a_{3 1}+a_{2
2}\,a_{4 0}\,a_{1 0}
 +a_{1 1}\,a_{4 1}\,a_{2 0}
  -a_{1 1}\,a_{3 2}\,a_{2 0} -a_{1 1}\,a_{3 0}\,a_{3 1}\\
 && +a_{1 1}\,a_{3 0}\,a_{2 2}-a_{2 0}\,a_{3 0}\,a_{2 2}
 -a_{1 0}\,a_{3 1}\,a_{2 2}-a_{1 0}\,a_{3 2}\,a_{3 0}
 -a_{2 1}\,a_{4 1}\,a_{1 0} -a_{4 1}\,a_{2 0}^2\\
&&{}
 +a_{1 0}\,a_{3 1}^2+a_{3 2}\,a_{2 0}^2
 +a_{3 0}\,a_{4 1}\,a_{1 0}+a_{2 0}\,a_{3 0}\,a_{3 1}, \\
&& {\bf a_{0 0}\,a_{4 4}}-a_{4 1}\,a_{3 0}+a_{4 0}\,a_{3 1}
 -a_{4 0}^2, \\
&& {\bf a_{0 0}\,a_{4 3},}-a_{4 1}\,a_{2 0}-a_{4 0}\,a_{3 0}
 +a_{4 0}\,a_{2 1}, \\
&& {\bf a_{0 0}\,a_{4 2}}-a_{4 0}\,a_{2 0}-a_{4 1}\,a_{1 0}
 +a_{4 0}\,a_{1 1}, \\
&& {\bf a_{0 0}\,a_{3 3}}-a_{4 1}\,a_{1 0}-a_{3 0}^2
 +a_{3 0}\,a_{2 1}+a_{4 0}\,a_{1 1}-a_{3 1}\,a_{2 0}, \\
&& {\bf a_{0 0}\,a_{3 2}}-a_{3 1}\,a_{1 0}+a_{4 0}\,a_{1 0}-a_{0
0}\,a_{4 1}
 -a_{3 0}\,a_{2 0}+a_{3 0}\,a_{1 1}, \\
&& {\bf a_{0 0}\,a_{2 2}}-a_{1 0}\,a_{2 1}+a_{1 0}\,a_{3 0}-a_{2
0}^2
 +a_{2 0}\,a_{1 1}-a_{3 1}\,a_{0 0} \rangle.
\end{array}
$$

We do not write the characteristic ideal
$\mathcal{I}_{\overline{a}}(K_6)$, because it contains 48
polynomials using the reduced polynomial with $a_{54}=0$. The
Gr\"obner basis of $\mathcal{I}(K_6)$ contains 104 polynomials.
The following are examples of polynomials $\Phi(x,y)$ such that
$K_n$ is strongly $\Phi$-polynomial.\\

\centerline{
\begin{tabular}{c|l}
Graph & Polynomial \\
\hline
 $K_3$ & $x^2 y^2+x^2+y^2+3 x+3 y+1$ \\
 $K_4$ & $x^3 y+x y^3+x^2 y^2+1 $ \\
 $K_5$ & $2 x^4 y^4+2 x^4+2 y^4+x^3 y+x y^3+x^2 y^2+1 $\\
 $K_6$ & $x^5 y^5+x^5+y^5-x^4 y^2-x^2 y^4-x^3 y^3+x+y+1$.\\
 \end{tabular}
}

To finish this section, two remarks. First, note that the complete
graph $K_{d+1}$ has $K_d$ as induced subgraph. Thus, a polynomial
with coefficients in
$\mathbb{V}(\mathcal{I}_{\overline{a}}(K_{d+1}))$ having all the
coefficients with some subindex $n$ equal to zero must be in
$\mathbb{V}(\mathcal{I}_{\overline{a}}(K_{d}))$. In terms of
ideals,
$$
\mathcal{I}_{\overline{a}}(K_{d+1})(a_{00},a_{10},\dots,a_{d-1,d-1},0,\dots,0)
=\mathcal{I}_{\overline{a}}(K_d).
$$
Second. As a consequence of Proposition~\ref{HUns}, if $K_d$ is
strongly $\Phi$-polynomial, then the polynomial $xy\,\Phi(x,y)$
satisfies the conditions of $K_{d+1}$. Therefore if we substitute
$a_{ij}$ by $a_{i+1\,j+1}$ in $\mathcal{I}_{\overline{a}}(K_d)$
the resulting ideal is contained in
$\mathcal{I}_{\overline{a}}(K_{d+1})$.

The above relations between the ideals $K_3$, $K_4$, $K_5$ and
$K_6$ can be checked.

\section{Concluding remarks}

We have presented and solved a number of questions
concerning polynomial graphs. Nevertheless, we are aware that there are
many open questions. Let us remark at least three of them.

The first one is, obviously, to prove or disprove the
  conjecture: Either to prove that if $(H,\Phi)$ is a standard pairing
  and $H$ is $\Phi$-polynomial, then $H$ is strongly
  $\Phi$-polynomial or to find a standard pairing $(H,\Phi)$ such that
  $H$ is $\Phi$-polynomial but not strongly $\Phi$-polynomial.

Second. We have seen that any strongly $\Phi$-polynomial graph is vertex
  transitive.  But all examples we have are Cayley graphs. Therefore,
  it is a natural question to ask if every strongly $\Phi$-polynomial
  graph is a Cayley graph.  In particular it would be interesting to
  know if Petersen's graph, which is vertex transitive but is not a Cayley
  graph, is polynomial (our guess is that it is not).

Third. The discussions in this paper are depending on the finiteness
  of $H$. It would be interesting to develop methods for $d$-regular graphs
  not necessarily finite, and generalize the conjecture.

\section{Acknowledgements}
We would like to thank Pelegr\'{\i} Viader for his many helpful
comments and his insightful perusal of our first draft.

\bibliographystyle{abbrv}

\end{document}